\documentclass{article}

\usepackage{amssymb,amsmath,bm,hyperref}
\usepackage{amssymb,amsmath,bm,geometry,graphics,graphicx,url,color}
\def\no{\noindent}
\def\pmatrix{\left(\begin{array}}
\def\endpmatrix{\end{array}\right)}
\def\udots{\reflectbox{$\ddots$}}  

\def\CC{\mathbb{C}}

\def\RR{\mathbb{R}}

\def\B{{\cal B}}
\def\C{{\cal C}}

\def\I{{\cal I}}

\def\P{{\cal P}}

\def\dd{\mathrm{d}}

\newtheorem{theo}{Theorem}
\newtheorem{lem}{Lemma}

\newtheorem{rem}{Remark}

\def\proof{\noindent\underline{Proof}\quad}
\def\QED{\mbox{~$\Box{~}$}}

\def\bfb{{\bm{b}}}
\def\bfc{{\bm{c}}}

\def\bfe{{\bm{e}}}

\def\bfgamma{{\bm{\gamma}}}

\begin{document}

\title{(Spectral) Chebyshev collocation methods for solving differential equations}

\author{Pierluigi\,Amodio\,\thanks{Dipartimento di Matematica, Universit\`a di Bari, Italy,
{\tt pierluigi.amodio@uniba.it}} \and
Luigi\,Brugnano\,\thanks{Dipartimento di Matematica e Informatica ``U.\,Dini'', Universit\`a di Firenze, Italy, {\tt luigi.brugnano@unifi.it}} \and
Felice\,Iavernaro\,\thanks{Dipartimento di Matematica, Universit\`a di Bari, Italy, {\tt felice.iavernaro@uniba.it}}}

\maketitle

\begin{abstract} Recently, the efficient numerical solution of Hamiltonian problems has been tackled by defining the class of energy-conserving Runge-Kutta methods named {\em Hamiltonian Boundary Value Methods (HBVMs)}. Their derivation relies on the expansion of the vector field along the Legendre orthonormal basis. Interestingly, this approach can be extended to cope with other orthonormal bases and, in particular, we here consider the case of the Chebyshev polynomial basis. The corresponding Runge-Kutta methods were previously obtained by Costabile and Napoli \cite{CoNa2001}. In this paper, the use of a different framework allows us to carry out a novel analysis of the methods also when they are used as spectral formulae in time, along with some generalizations of the methods.

\bigskip
\no{\bf Keywords:} Legendre polynomials, Hamiltonian Boundary Value Methods, HBVMs, Chebyshev polynomials, Chebyshev collocation methods.

\bigskip
\no{\bf AMS-MSC:} 65L06, 65L05.

\end{abstract}

\section{Introduction}\label{sec1}

The numerical solution of Hamiltonian problems has been the subject of many investigations in the last decades, due to the fact that Hamiltonian problems are not structurally stable against generic perturbations, like those induced by a general-purpose numerical method  used to approximate their solutions. 
The main features of a canonical Hamiltonian system are surely the symplecticness of the map and the conservation of energy, which cannot be simultaneously maintained by any given Runge-Kutta or B-series method \cite{CFM2006}. Consequently, numerical methods have been devised in order of either  defining a symplectic discrete map, giving rise to the class of {\em symplectic methods} (see, e.g., the monographs \cite{SSC1994,LeRe2004,GNI2006,BlCa2016}, and references therein, the review paper \cite{SS2016}, and related approaches \cite{But21}), or being able to conserve the energy, resulting in the class of {\em energy-conserving methods} (see, e.g., \cite{LQR1999,IaPa2007,IaPa2008,IaTr2009,CMcLMcLOQW2009,BIT2009,BIT2010,H2010,BIS2010,LW2016} and the monographs \cite{LIMBook2016}).   In particular, {\em Hamiltonian Boundary Value Methods (HBVMs)} \cite{BIT2009,BIT2010,BIS2010,LIMBook2016,BIT2012} (as well as their generalizations, see, e.g. \cite{BCMR2012,BI2012,BIT2012_2,BS2014,ABI2015,BGIW2018,BMR2019,BIZ2020,ABI2022}) form a special class of  Runge-Kutta methods with a (generally)  rank-deficient coefficient matrix.  HBVMs, in turn, can be also viewed as the outcome obtained after a local projection of the vector field onto a finite-dimensional function space: in particular, the set of polynomials of a given degree. For this purpose, the Legendre orthonormal polynomial basis has been considered so far \cite{BIT2012_1} (see also \cite{ABI2019_0}), and their use as spectral methods in time has been also investigated both theoretically and numerically \cite{BMR2019,BIMR2019,BGZ2019,ABI2020,ABIM2020,BBTZ2020}.
Remarkably, as was already observed in \cite{BIT2012_1}, this idea is even more general and can be adapted to other finite-dimensional function spaces and/or different bases. Following this route, in this paper we consider a class of Runge-Kutta methods based on the use of the Chebyshev polynomial basis. 

It turns out that, with this choice, our approach finally leads us back to the same formulae introduced by Costabile and Napoli  by considering the classical collocation conditions based on the Chebyshev abscissae \cite{CoNa2001}. Exploiting the Routh-Hurwitz criterion, in \cite{CoNa2004} they derived the A-stability property of the formule up to order 20  while,  more recently, they also extended the methods to cope with  $k$-th order  problems \cite{CoNa2011}. 

In this context, we are mainly interested in a spectral-time implementation of these formulae, which means that the given approximation accuracy (usually close to the machine epsilon)  is achieved by increasing the order of the formulae rather than reducing the integration stepsize. An interesting consequence of such a strategy is that, modulo the effects of round-off errors, the numerical solution will mimic the exact solution so that, when applied to a Hamiltonian system the method will inherit both the simplecticity and the conservation properties.    

The advantage of using the Chebyshev basis stems from the fact that all the entries in the Butcher tableau of the corresponding Runge-Kutta methods can be given in closed form, thus avoiding the introduction of round-off errors when numerically computing them (as is the case with the Legendre basis, where the Gauss-Legendre nodes need to be numerically computed). In this respect, not only does the analysis provided within the new framework help in deriving the explicit expression of the coefficients of the methods for any arbitrary high order, but it is also very useful to discuss the convergence rate when they are used as spectral methods in time and to derive an efficient implementation strategy based on the discrete cosine transform: these latter aspects are new, at the best of our knowledge. Further, also a generalization of the methods is sketched, similar to that characterizing HBVMs.

With this premise, the structure of the paper is as follows: in Section~\ref{Hamil} we provide a brief introduction to the framework used to derive HBVMs, relying on the use of the Legendre polynomial basis;  in Section~\ref{ceby} the approach is extended to cope with the Chebyshev polynomial basis; in Section~\ref{anal} a thorough analysis of the resulting methods is given, also when they are used as spectral methods in time; in Section~\ref{numtest} we provide some numerical tests, to assess the theoretical findings; at last a few conclusions are given in Section~\ref{fine}.  

\section{Hamiltonian Boundary Value Methods (HBVMs)}\label{Hamil}

In order to introduce HBVMs, let us consider a canonical Hamiltonian problem in the form
\begin{equation}\label{Hprob}
\dot y = J \nabla H(y), \qquad y(0) = y_0\in\RR^{2m},  \qquad J=\pmatrix{cc} &I_m\\-I_m\endpmatrix = -J^\top,
\end{equation}
with $H:\RR^{2m}\rightarrow\RR$ the Hamiltonian function (or {\em energy}) of the system.\footnote{In fact, for isolated mechanical systems, $H$ has the physical meaning of the total energy.} As is easily understood $H$ is a constant of motion, since
$$\frac{\dd}{\dd t}H(y) = \nabla H(y)^\top\dot y =  \nabla H(y)^\top J \nabla H(y)=0,$$
being $J$ skew-symmetric. The simple idea, on which HBVMs rely on, is that of reformulating the previous conservation property in  terms of the line integral of $\nabla H$ along the path defined by the solution $y(t)$:
$$H(y(t)) = H(y_0) + \int_0^t \nabla H(y(\tau))^\top\dot y(\tau)\dd\tau.$$
Clearly, the integral at the right-hand side of the previous equality vanishes, because the integrand is identically zero by virtue of (\ref{Hprob}), so that the conservation holds true for all $t>0$. The solution of (\ref{Hprob}) is the unique function satisfying such a conservation property for all $t>0$. However, if we consider a discrete-time dynamics, ruled by a {\em timestep} $h$, there exist infinitely many paths $\sigma$ such that:
\begin{eqnarray}\nonumber
\sigma(0) &=& y_0, \qquad \sigma(h) =: y_1\approx y(h),\\
0&=&h\int_0^1 \nabla H(\sigma(ch))^\top\dot\sigma(ch)\dd c.\label{lim}
\end{eqnarray}
The path $\sigma$ obviously defines a one-step numerical method that conserves the energy, since 
$$H(y_1)~=~H(y_0) + h\int_0^1 \nabla H(\sigma(ch))^\top\dot\sigma(ch)\dd c ~=~ H(y_0),$$
even though now the integrand is no more identically zero. The methods derived in this framework have been named {\em line integral methods}, due to the fact that the path $\sigma$ is defined so that the line integral in (\ref{lim}) vanishes. Line integral methods have been thoroughly analyzed in the monograph \cite{LIMBook2016} (see also the review paper \cite{BI2018}). Clearly, in the practical implementation of the methods, the integral is replaced by a quadrature of enough high-order, thus providing a fully discrete method, even though we shall not consider, for the moment, this aspect, for which the reader may refer to the above mentioned references.\footnote{This amounts to study HBVMs as {\em continuous-stage Runge-Kutta methods}, as it has been done in \cite{ABI2019_0,ABI2022_0}.}

Interestingly enough, after some initial attempts to derive methods in this class \cite{IaPa2007,IaPa2008,IaTr2009,BIT2009,BIT2015}, a systematic way for their derivation was found in \cite{BIT2012_1}, which is based on a local Fourier expansion of the vector field in (\ref{Hprob}). In fact, by setting 
\begin{equation}\label{fy}
f(y)=J\nabla H(y), 
\end{equation}
and using hereafter, depending on the needs, either one or the other notation, one may rewrite problem (\ref{Hprob}), on the interval $[0,h]$, as:
\begin{equation}\label{y1}
\dot y(ch) = \sum_{j\ge0} P_j(c)\gamma_j(y), \qquad c\in[0,1], \qquad y(0) = y_0, 
\end{equation}
where $\{P_j\}_{j\ge0}$ is the Legendre orthonormal polynomial basis on $[0,1]$, 
\begin{equation}\label{leg}
P_i\in\Pi_i, \qquad \int_0^1 P_i(\tau)P_j(\tau)\dd \tau = \delta_{ij}, \qquad i,j\ge0,
\end{equation}
$\Pi_i$, hereafter, denotes the vector space of polynomials of degree at most $i$, $\delta_{ij}$ is the Kronecker symbol, and 
\begin{equation}\label{gamj}
\gamma_j(y) = \int_0^1P_j(\tau)f(y(\tau h))\dd\tau, \qquad j \ge 0,
\end{equation} 
are the corresponding Fourier coefficients. 
The solution of the problem is formally obtained,  in terms of the unknown Fourier coefficients, by integrating both sides of Equation~(\ref{y1}):
\begin{equation}\label{y}
y(ch) = y_0 + h\sum_{j\ge0}\int_0^c P_j(x)\dd x\, \gamma_j(y), \qquad c\in[0,1].
\end{equation}
A polynomial approximation $\sigma\in\Pi_s$ can be formally obtained by truncating the previous series to finite sums:
\begin{equation}\label{sig1}
\dot\sigma(ch) = \sum_{j=0}^{s-1} P_j(c) \gamma_j(\sigma), \qquad c\in[0,1], \qquad \sigma(0)=y_0,
\end{equation}
and
\begin{equation}\label{sig}
\sigma(ch) = y_0 + h\sum_{j=0}^{s-1}\int_0^c P_j(x)\dd x\, \gamma_j(\sigma), \qquad c\in[0,1],
\end{equation}
respectively, with $\gamma_j(\sigma)$ defined according to (\ref{gamj}), upon replacing $y$ with $\sigma$. Whichever the degree $s\ge1$ of the polynomial approximation, the following result holds true, where the same notation used above holds.

\begin{theo}
\label{Hcons} 
$H(y_1)=H(y_0)$. 
\end{theo}
\proof In fact, one has:
\begin{eqnarray*} 
\lefteqn{H(y_1)-H(y_0) ~=~ H(\sigma(h))-H(\sigma(0)) ~=~ h\int_0^1 \nabla H(\sigma(ch))^\top\dot\sigma(ch)\dd c}\\
&=& h\int_0^1 \nabla H(\sigma(ch))^\top\sum_{j=0}^{s-1} P_j(c)\gamma_j(\sigma)\dd c~
=~ h\sum_{j=0}^{s-1} \left[ \int_0^1 P_j(c)\nabla H(\sigma(ch))\dd c\right]^\top \gamma_j(\sigma)\\
&=& h\sum_{j=0}^{s-1} \left[ \int_0^1 P_j(c)\nabla H(\sigma(ch))\dd c\right]^\top J  \left[ \int_0^1 P_j(c)\nabla H(\sigma(ch))\dd c\right] ~=~ 0,
\end{eqnarray*}
due to the fact that $J$ is skew-symmetric.\,\QED\bigskip

The next result states that the method has order $2s$.

\begin{theo}\label{ord2s}
$y_1-y(h)=O(h^{2s+1})$. 
\end{theo}
\proof See \cite[Theorem~1]{BIT2012_1}.\,\QED\bigskip

It must be emphasized that, when the method is used as a spectral method in time, then the concept of order does not hold anymore. Instead, the following result can be proved, under regularity assumptions on $f$ in (\ref{fy}).

\begin{theo}\label{spectral}
Let us assume $f(\sigma(t))$ be analytical in a closed ball of radius $r^*$ centered at $0$, then  for all $h\in(0,h^*]$, with $h^*<r^*$, there exist $M=M(h^*)>0$ and $\rho>1$, $\rho\sim h^{-1}$, such that:\,\footnote{Hereafter, $|\cdot|$ will devote any convenient vector norm.}
\begin{equation}\label{norma}
\|\sigma-y\|~:=~\max_{c\in[0,1]}|\sigma(ch)-y(ch)| ~\le~ h M\rho^{-s}.
\end{equation}
\end{theo}
\proof See \cite[Theorem~2]{ABI2020}.\,\QED \bigskip

\begin{rem}\label{serve tutto}
When using HBVMs as spectral methods in time, it is more important considering the measure of the error (\ref{norma}), rather than the error at $h$. In fact, the timestep used, in such a case, may be large, and even huge, so that the whole approximation in the interval $(0,h)$ is needed. 
\end{rem}

\begin{rem}\label{formachiusa}
As previously stated, an important issue,  when using numerical methods as spectral methods, stems from the fact that one needs the relevant coefficients be computed for very high-order formulae. When some of such coefficients are evaluated numerically, as is the case for the abscissae of the Gauss-Legendre formulae, this may introduce errors that, even though small, may affect the accuracy of the resulting method. It is then useful to have formulae for which all the involved coefficients are known in closed form, and this motivates the present paper.
\end{rem}

\section{Chebyshev-Runge-Kutta methods}\label{ceby}

An interesting way to interpret the approximation (\ref{sig1}) is that of looking for the coefficients 
$\gamma_0,\dots,\gamma_{s-1}$, in the polynomial approximation 
\begin{equation}\label{sig1_new}
\dot\sigma(ch) = \sum_{j=0}^{s-1} P_j(c)\gamma_j, \qquad c\in[0,1],\qquad \sigma(0) = y_0\in\RR^m,
\end{equation}
to\,\footnote{For sake of simplicity, hereafter we shall assume $f$ be an analytic function.} 
\begin{equation}\label{ode1}
\dot y(ch) = f(y(ch)), \qquad c\in[0,1], \qquad y(0)=y_0\in\RR^m,
\end{equation} 
such that the residual function
\begin{equation}\label{rch}
r(ch) ~:=~ \dot\sigma(ch) - f(\sigma(ch)), \qquad c\in[0,1],
\end{equation}
be orthogonal to $\Pi_{s-1}$. That is,
\begin{equation}\label{orto1}
0 ~=~ \int_0^1 P_i(c)r(ch)\dd c ~=~ \int_0^1 P_i(c)\left( \dot\sigma(ch) - f(\sigma(ch))\right)\dd c, \qquad i=0,\dots,s-1.
\end{equation}
By considering the orthogonality relations (\ref{leg}), this amounts to require that
\begin{equation}\label{gami0}
\gamma_i = \int_0^1 P_i(c) f\left( y_0+h\sum_{j=0}^{s-1}\int_0^cP_i(x)\dd x\gamma_j\right)\dd c, \qquad i=0,\dots,s-1,\end{equation}
namely (\ref{gamj}) with $\sigma$ replacing $y$, and the new approximation given by 
\begin{equation}\label{ynew}
y_1~:=~\sigma(h) ~\equiv~ y_0 + h\sum_{j=0}^{s-1}\int_0^1 P_j(x)\dd x\gamma_j.
\end{equation}
Approximating the integrals in (\ref{gami0}) by a Gauss-Legendre quadrature of order $2k$, with $k\ge s$, then provides a HBVM$(k,s)$ method \cite{LIMBook2016,BI2018,BIT2010,BIT2012_1}.

A generalization of the orthogonality requirement (\ref{gami0}) consists in considering a suitable weighting function
\begin{equation}\label{wc}
\omega(c)\ge 0, \qquad c\in[0,1], \qquad \int_0^1 \omega(c)\dd c=1,
\end{equation}
and a polynomial basis orthonormal w.r.t. the induced product,
\begin{equation}\label{legw}
P_i\in\Pi_i, \qquad \int_0^1 \omega(\tau)P_i(\tau)P_j(\tau)\dd \tau = \delta_{ij}, \qquad i,j\ge0,
\end{equation}
then requiring
\begin{equation}\label{ortow}
0 ~=~ \int_0^1 \omega(c) P_i(c)r(ch)\dd c ~\equiv~ \int_0^1 \omega(c)P_i(c)\left( \dot\sigma(ch) - f(\sigma(ch))\right)\dd c, \qquad i=0,\dots,s-1.
\end{equation}
Consequently, now the coefficients of the polynomial approximation (\ref{sig1_new}) turn out to solve the following set of equations,
\begin{equation}\label{gamiw}
\gamma_i = \int_0^1 \omega(c)P_i(c) f\left( y_0+h\sum_{j=0}^{s-1}\int_0^cP_i(x)\dd x\gamma_j\right)\dd c, \qquad i=0,\dots,s-1,\end{equation}
in place of (\ref{gami0}).\footnote{Clearly,  (\ref{gami0}) is derived by considering the special case $\omega(c)\equiv1$.}
In so doing, we obtain a polynomial approximation formally still given by (\ref{sig}), with the Fourier coefficients now defined as
\begin{equation}\label{gamiws}
\gamma_j ~\equiv~\gamma_i(\sigma) ~:=~ \int_0^1 \omega(c)P_i(c) f(\sigma(ch))\dd c, \qquad i=0,\dots,s-1.
\end{equation}
Hereafter, we shall consider the weighting function
\begin{equation}\label{cebyw}
\omega(c) = \frac{1}{\pi \sqrt{c(1-c)}}, \qquad c\in(0,1),
\end{equation}
which satisfies (\ref{wc}) and provides the shifted and scaled Chebyshev polynomials of the first kind, i.e.,\footnote{As is usual, hereafter, $T_j(x)$, $j\ge0$, denote the Chebyshev polynomials of the first kind.}
\begin{equation}\label{cebypol}
P_0(c) ~=~ T_0(2c-1) ~\equiv~ 1, \qquad P_j(c) ~=~ \sqrt{2}\,T_j(2c-1), \quad j\ge1, \qquad c\in[0,1],
\end{equation}
satisfying (\ref{legw}).

For later use, we also report the relations between such polynomials and their integrals:
\begin{eqnarray}\nonumber
\int_0^c P_0(x)\dd x &=& \frac{1}2\left[ \frac{P_1(c)}{\sqrt{2}}+P_0(c) \right],\\ \label{cebyint}
\int_0^c P_1(x)\dd x &=& \frac{1}8\left[  P_2(c)-\sqrt{2}P_0(c) \right],\\ \nonumber
\int_0^c P_j(x)\dd x &=&  \frac{1}4\left[  \frac{P_{j+1}(c)}{j+1} -   \frac{P_{j-1}(c)}{j-1}  - \frac{(-1)^j2\sqrt{2}P_0(c)}{j^2-1}\right] , \qquad j\ge 2.
\end{eqnarray}
In particular, for $c=1$ one obtains, by considering that $P_0(1) = 1$ and $P_j(1) = \sqrt{2}$, for all $j\ge1$:
\begin{equation}\label{intP1}
\int_0^1 P_j(x)\dd x = \left\{
\begin{array}{cl}
1, &~j~=~0,\\[2mm]
0, &~j\quad \mbox{odd,}\\[2mm]
\frac{\sqrt{2}}{1-j^2}, &~j\quad\mbox{even}, \quad j\ge2.
\end{array}\right.
\end{equation}

\subsection{Discretization}

As is clear, the integrals involved in (\ref{gamiw})-(\ref{gamiws}) cannot be computed exactly, and need to be approximated by using a numerical method: this latter can be naturally chosen as the Gauss-Chebysvev interpolatory quadrature formula of order $2s$ on the interval $[0,1]$, with nodes
\begin{equation}\label{ci}
c_i = \frac{1+\cos\theta_i}2, \qquad \theta_i = \frac{2i-1}{2s}\pi, \qquad i=1,\dots,s,
\end{equation}
and weights $$\omega_i = \frac{1}s, \qquad i=1,\dots,s.$$ 
We recall that $P_s(c_i)=0$, $i=1,\dots,s$, due to fact that (see (\ref{cebypol})), $x_i=\cos\theta_i$, $i=1,\dots,s$ are the roots of $T_s(x)$. Consequently, the Fourier coefficients (\ref{gamiw})-(\ref{gamiws}) will be approximated as:
\begin{equation}\label{gamis}
\hat\gamma_i = \frac{1}s\sum_{j=1}^s P_i(c_j)f\left( y_0+h\sum_{\ell=0}^{s-1}\int_0^{c_j}P_\ell(x)\dd x\,\hat\gamma_\ell\right)\dd c, \qquad i=0,\dots,s-1.\end{equation}
In so doing, we obtain a new polynomial approximation, 
\begin{equation}\label{uch}
u(ch) = y_0 + h\sum_{j=0}^{s-1}\int_0^cP_j(x)\dd x\,\hat\gamma_j, \qquad c\in[0,1],
\end{equation}
in place of $\sigma$. 
Setting $Y_j$ the argument of $f$ in Equation (\ref{gamis}), and taking into account (\ref{uch}), we arrive at the $s$-stage Runge-Kutta method with stages
\begin{equation}\label{Yi}
Y_i ~:=~u(c_ih) ~=~ y_0 + h\sum_{j=1}^s\underbrace{\left(\frac{1}s\sum_{\ell=0}^{s-1} \int_0^{c_i} P_\ell(x)\dd x P_\ell(c_j)\right)}_{=:\,a_{ij}} f(Y_j), \qquad i=1,\dots,s, 
\end{equation}
with the new approximation given by:
\begin{equation}\label{y1dig}
y_1 ~:=~ u(h) = y_0 + h\sum_{i=1}^s \underbrace{\left(\frac{1}s\sum_{\ell=0}^{s-1} \int_0^1 P_\ell(x)\dd x P_\ell(c_i)\right)}_{=:\,b_i} f(Y_i).
\end{equation}
Consequently, the abscissae and weights $(c_i,b_i)$, $i=1,\dots,s$, and the entries $a_{ij}$, $i,j=1,\dots,s$, define an $s$-stage Runge-Kutta method. An explicit expression of the abscissae has been given in (\ref{ci}). Let us now derive more refined expressions for the weights and the Butcher matrix. For this purpose, let us define the following matrices:
\begin{equation}\label{PIO}
\P_s = \Big( P_{j-1}(c_i)\Big)_{i,j=1,\dots,s}, \qquad \I_s = \left( \int_0^{c_i} P_{j-1}(x)\dd x\right)_{i,j=1,\dots,s}, \qquad \Omega = \frac{1}s I_s,
\quad \in~\RR^{s\times s},
\end{equation}
with $I_s$ the identity matrix, and
\begin{equation}\label{Xs}
X_s = \pmatrix{ccccccc}
\frac{1}2              & -\sqrt{2}\beta_2    & \alpha_3 & \alpha_4 &\dots  & \alpha_{s-1} &\alpha_s\\[1mm] 
\sqrt{2}\beta_1      &   0                       &  -\beta_1   & \\
                            & \beta_2             &  0        & -\beta_2 & &O\\
                            &                            & \beta_3    &0 &-\beta_3\\
                            &                           &         &\ddots  &\ddots &\ddots\\
                            &O                          &           &   &\beta_{s-2} &0 &-\beta_{s-2}\\
                            &                            &          &          &  &\beta_{s-1} & 0
\endpmatrix        \in\RR^{s\times s},                   
\end{equation}
with
\begin{equation}\label{Xs1}
\beta_j = \frac{1}{4j}, \quad j\ge 1,\qquad  \alpha_j = (-1)^j8\sqrt{2}\beta_j\beta_{j-2}, \quad j\ge 3.  
\end{equation}

The following results hold true.

\begin{lem}\label{PIOX}
With reference to (\ref{ci}) and the matrices (\ref{PIO}) and (\ref{Xs}), one has:
$$ (\P_s)_{ij} =  \left[\sqrt{2} + \delta_{j1}(1-\sqrt{2})\right] \cos (j-1)\theta_i, \qquad i,j=1,\dots,s,$$
$$\I_s = \P_s X_s, \qquad P_s^\top\Omega \P_s = I_s.$$
\end{lem}
\proof
The first statement follows from (\ref{cebypol}) and the well-known fact that 
$$T_{j-1}(\cos\theta_i) = \cos(j-1)\theta_i.$$
The second statement follows from (\ref{cebyint}), by considering that $P_s(c_i)=0$, $i=1,\dots,s$. Finally, one has,
by setting, hereafter, $e_i\in\RR^s$ the $i$th unit vector:
\begin{eqnarray*}
e_i^\top\left( \P_s^\top\Omega \P_s\right) e_j &=& \frac{1}s(\P_s e_i)^\top(\P_se_j) \\[1mm]
&=& \frac{\left[\sqrt{2} + \delta_{i1}(1-\sqrt{2})\right]\left[\sqrt{2} + \delta_{j1}(1-\sqrt{2})\right]}s
\sum_{k=1}^s\cos(i-1)\theta_k \cdot \cos(j-1)\theta_k \\
&=& \delta_{ij},
\end{eqnarray*}
due to the fact that, for all $i,j=1,\dots,s$,
$$\sum_{k=1}^s\cos(i-1)\theta_k \cdot \cos(j-1)\theta_k ~=~ \left\{
\begin{array}{cl}
0, & i\ne j,\\[1mm]
s, & i=j=1,\\[1mm]
\frac{s}2, &i=j\ge 1.\,\QED
\end{array}\right.$$\smallskip

Consequently, we can now state the following result.

\begin{theo}\label{buttabc}
The Butcher matrix of the Runge-Kutta method (\ref{Yi})-(\ref{y1dig}) is given by
\begin{equation}\label{butA}
A ~\equiv~(a_{ij})~ :=~ \I_s\P_s^\top\Omega ~\equiv~ \P_s X_s \P_s^\top\Omega ~\equiv~  \P_s X_s \P_s^{-1}.
\end{equation}
The corresponding weights are given by:
\begin{equation}\label{bi}
b_i := \frac{1}s\left[ 1 - 2\sum_{j=1}^{\lceil s/2\rceil-1} \frac{ \cos \frac{2i-1}s j\pi}{4j^2-1}\right], \qquad i=1,\dots,s.
\end{equation}
\end{theo}
\proof
From Lemma~\ref{PIOX}, the last two equalities in (\ref{butA}) easily follow. The first equality in (\ref{butA}) then follows from (\ref{Yi}), by  taking into account (\ref{PIO}):
$$a_{ij}~=~e_i^\top A e_j ~=~ e_i^\top \I_s\P_s^\top \Omega e_j ~=~ \frac{1}s (e_i^\top \I_s)(e_j^\top \P_s)^\top
~=~ \frac{1}s\sum_{\ell=0}^{s-1} \int_0^{c_i}P_\ell(x)\dd x P_\ell(c_j),$$
which coincides with the formula given in  (\ref{Yi}).
The formula (\ref{bi}) of the weights follows by considering that from (\ref{y1dig}), taking into account (\ref{intP1}) and (\ref{ci}), one has:

\begin{eqnarray*}
b_i &=&\frac{1}s\sum_{\ell=0}^{s-1} \int_0^1 P_\ell(x)\dd x P_\ell(c_i) ~=~ \left(\int_0^1P_0(c)\dd c\,,\dots,\int_0^1P_{s-1}(c)\dd c\right)^\top \P_s^\top\Omega e_i\\
&=& \frac{1}s\left( 1,\,0,\,\frac{-\sqrt{2}}{2^3-1},\,0,\,\frac{-\sqrt{2}}{4^3-1},\,0,\,\frac{-\sqrt{2}}{6^3-1},\,\dots\right)
\pmatrix{c} 1\\ \sqrt{2}\cos\theta_i \\ \vdots \\ \sqrt{2}\cos(s-1)\theta_i\endpmatrix\\[1mm]
&=& \frac{1}s\left[ 1 - 2\sum_{j=1}^{\lceil s/2\rceil-1} \frac{ \cos 2j\theta_i}{4j^2-1}\right]~=~\frac{1}s\left[ 1 - 2\sum_{j=1}^{\lceil s/2\rceil-1} \frac{ \cos \frac{2i-1}s j\pi}{4j^2-1}\right].\,\QED
\end{eqnarray*}\smallskip 

Moreover, concerning the weights of the method, the following result holds true.

\begin{theo}\label{bipos} The weights $\{b_i\}_{i=1,\dots,s}$ defined in (\ref{bi}) are all positive.\end{theo}  
\proof In fact, from (\ref{bi}), having fixed a given value of $s\ge 1$, one has:
\begin{eqnarray*}
b_i &=&\frac{1}s\left[ 1 - 2\sum_{j=1}^{\lceil s/2\rceil-1} \frac{ \cos \frac{2i-1}s j\pi}{4j^2-1}\right]
~\ge~ \frac{1}s\left[ 1 - \sum_{j=1}^{\lceil s/2\rceil-1} \frac{2}{4j^2-1}\right]\\
&=& \frac{1}s\left[ 1 - \sum_{j=1}^{\lceil s/2\rceil-1} \left(\frac{1}{2j-1}-\frac{1}{2j+1}\right)\right]
~=~\frac{1}s\left[ 1 -1 + \frac{1}{2(\lceil s/2\rceil-1)+1}\right]\\
&=&  \frac{1}s\,\frac{1}{2\lceil s/2\rceil-1} ~\ge~\frac{1}s\,\frac{1}s ~=~\frac{1}{s^2}.\,\QED
\end{eqnarray*}\smallskip

In conclusion, we have derived the family of $s$-stage Runge-Kutta method whose Butcher tableau, with reference to (\ref{ci}), (\ref{bi}), and (\ref{PIO})--(\ref{Xs1}), is given by:

\begin{equation}\label{tabCCM}
\begin{array}{c|c}
\bfc & \I_s \P_s^\top\Omega\\ \hline \\[-3mm]
& \bfb^\top
\end{array} 
\quad \equiv \quad \begin{array}{c|c}
\bfc & \P_s X_s\P_s^\top\Omega\\ \hline \\[-3mm]
& \bfb^\top
\end{array}
\quad \equiv \quad \begin{array}{c|c}
\bfc & \P_s X_s\P_s^{-1}\\ \hline \\[-3mm]
& \bfb^\top
\end{array}
\end{equation}
having set, as is usual,
\begin{equation}\label{bc}
\bfc=(c_1,\,\dots,\,c_s)^\top, \qquad \bfb=(b_1,\,\dots,\,b_s)^\top.
\end{equation}
In particular, the last form in (\ref{tabCCM}) can be regarded as a kind of $W$-transformation \cite{HW1996} of the method.

\begin{rem} The Runge-Kutta methods in (\ref{tabCCM})-(\ref{bc}) coincide with the {\em Chebyshev Collocation Methods} derived by Costabile and Napoli in \cite{CoNa2001}. For sake of brevity, we shall refer to the $s$-stage method as CCM$(s)$.\end{rem}

\subsection{Collocation methods}\label{collo}
As anticipated above, CCM$(s)$ methods are {\em collocation methods}. In fact, from (\ref{gamis}), (\ref{Yi}), and (\ref{tabCCM}), by setting
$$Y ~=~ \pmatrix{c} u(c_1h)\\ \vdots\\ u(c_sh)\endpmatrix ~\equiv~ u(\bfc h), \qquad \hat{\bfgamma}=\pmatrix{c}\hat\gamma_0\\ \vdots\\ \hat\gamma_{s-1}\endpmatrix,$$
the vector of the stages and of the Fourier coefficients, respectively,
one derives that:
$$\dot u(\bfc h) = \P_s\otimes I_m \hat{\bfgamma},\qquad  \hat{\bfgamma} = \P_s^\top\Omega\otimes I_m f(u(\bfc h)).$$
By combining the last two equations, and taking into account that, by virtue of Lemma~\ref{PIOX}, 
$\P_s\P_s^\top\Omega = I_s$, one eventually derives that
$$\dot u(\bfc h) = f(\bfc h),$$
i.e.,
$$\dot u(c_ih) = f(u(c_ih)), \qquad i=1,\dots,s.$$
Further, the collocation polynomial $u$ satisfies ~$u(0)=y_0$~ and ~$u(h)=:y_1$,~ as seen above.

\subsection{Interesting implementation details}\label{itera}

An interesting property of the Butcher matrix in the tableau (\ref{tabCCM}) derives from the fact that the multiplications
$$\sqrt{s}\, \P_s^\top \Omega \otimes I_m V \equiv \frac{1}{\sqrt{s}}\P_s^\top\otimes I_m V, \qquad \mbox{and} \qquad \frac{1}{\sqrt{s}}\P_s\otimes I_m Z,$$
with $V,Z\in\RR^{sm}$ given vectors, amount to the discrete cosine transform of $V$, ${\tt dct}(V)$, and the inverse discrete cosine transform of $Z$, ${\tt idct}(Z)$, respectively.\footnote{Here, we have used the name of the Matlab$^\copyright$  functions implementing the two transformations.} Consequently, a fixed-point iteration for computing the stages of the CCM$(s)$ method reads, by setting hereafter the vector ~$\bfe = \left( 1,\, 1,\,\dots,\,1\right)^\top\in\RR^s$:
\begin{equation}\label{dct}
Y^{(r+1)} = \bfe\otimes y_0 + h\,{\tt idct}\left( X_s\otimes I_m\,{\tt dct}(f(Y^{(r)}))\right), \qquad r=0,1,\dots,
\end{equation} 
which can be advantageous, for large values of $s$. In fact, by considering that the matrix $X_s$ in (\ref{Xs}) is {\em sparse}, the complexity of one iteration amounts to computing the vector field in the stage values, plus $O(ms\log(s))$ flops,\footnote{As is usual, 1 flop denotes a basic {\em fl}oating-point {\em op}eration.} whereas the standard implementation would require $2ms^2$ flops.

\subsection{Generalizations}\label{chbvmks}
We observe that a possible generalizations of the collocation methods described above is that of using, in the discretization procedure of the integrals involved in (\ref{gamiw})-(\ref{gamiws}), a Gauss-Chebysvev interpolatory quadrature formula of order $2k$ on the interval $[0,1]$, with nodes and weights
\begin{equation}\label{kges}
c_i = \frac{1+\cos\theta_i}2, \qquad \theta_i = \frac{2i-1}{2k}\pi, \qquad \omega_i = \frac{1}k, \qquad i=1,\dots,k,
\end{equation}
for a convenient $k\ge s$. In such a case, for $k>s$, the methods are no more collocation methods. This is analogous to what happens for HBVM$(k,s)$ methods \cite{LIMBook2016} that, when $k=s$, reduce to the $s$-stage Gauss-Legendre collocation methods but, for $k>s$, are no more a collocation methods. 
By choosing the absissae (\ref{kges}), the matrices (\ref{PIO})-(\ref{Xs})  respectively become:
\begin{equation}\label{PIO1}
\P_s = \Big( P_{j-1}(c_i)\Big)_{\scriptsize\begin{array}{l}i=1,\dots,k\\j=1,\dots,s\end{array}}, ~ 
\I_s = \left( \int_0^{c_i} P_{j-1}(x)\dd x\right)_{\scriptsize\begin{array}{l}i=1,\dots,k\\j=1,\dots,s\end{array}}\in\RR^{k\times s}, 
\qquad \Omega = \frac{1}k I_k \in\RR^{k\times k},
\end{equation}
and (see (\ref{Xs}))
$$
\hat X_s = \pmatrix{ccccccc}
\frac{1}2              & -\sqrt{2}\beta_2    & \alpha_3 & \alpha_4 &\dots  & \alpha_{s-1} &\alpha_s\\[1mm] 
\sqrt{2}\beta_1      &   0                       &  -\beta_1   & \\
                            & \beta_2             &  0        & -\beta_2 & &O\\
                            &                            & \beta_3    &0 &-\beta_3\\
                            &                           &         &\ddots  &\ddots &\ddots\\
                            &O                          &           &   &\beta_{s-2} &0 &-\beta_{s-2}\\
                            &                            &          &          &  &\beta_{s-1} & 0\\
                            \hline
                            &                            &          &          & & &\beta_s
\endpmatrix    \equiv \pmatrix{c} X_s \\ \hline 0,\dots,0,\beta_s\endpmatrix    \in\RR^{s+1\times s}.                   
$$ This latter matrix is such that, for $k>s$:
\begin{equation}\label{IsPs}
\I_s = \P_{s+1}\hat X_s,
\end{equation}
where $\P_{s+1}\in\RR^{k\times s+1}$ is defined similarly as in (\ref{PIO1}).\footnote{Clearly, when $k=s$, (\ref{IsPs}) reduces to 
$\I_s=\P_sX_s$, according to Lemma~\ref{PIOX}.} 
Finally, the corresponding Butcher tableau become (compare with (\ref{tabCCM})-(\ref{bc})):
$$
\begin{array}{c|c}
\bfc & \I_s \P_s^\top\Omega\\ \hline \\[-3mm]
& \bfb^\top
\end{array} 
\quad \equiv \quad \begin{array}{c|c}
\bfc & \P_{s+1} \hat X_s\P_s^\top\Omega\\ \hline \\[-3mm]
& \bfb^\top
\end{array}\, ,\qquad \bfc=(c_1,\,\dots,\,c_k)^\top, \qquad \bfb=(b_1,\,\dots,\,b_k)^\top,
$$
with the weights given by 
$$
b_i := \frac{1}k\left[ 1 - 2\sum_{j=1}^{\lceil s/2\rceil-1} \frac{ \cos \frac{2i-1}k j\pi}{4j^2-1}\right], \qquad i=1,\dots,k,
$$
in place of (\ref{bi}). Nevertheless, provided that $k\ge s$ holds, the analysis of such methods is similar to that of CCM$(s)$ methods (all the results in Section~\ref{anal} continue formally to hold), so that we shall not consider them further.

\section{Analysis of the methods}\label{anal}

In this section we carry out the analysis of the CCM$(s)$ method defined by the Butcher tableau (\ref{tabCCM})-(\ref{bc}). In particular, we study:
\begin{enumerate}
\item the symmetry of the method;
\item its order of convergence;
\item its linear stability;
\item its accuracy when used as a spectral method in time. 
\end{enumerate}
The analysis of items 2 and 3  uses different arguments from those provided in \cite{CoNa2004,CoNa2011}, whereas the analysis related to items 1 and 4 is novel.

\subsection{Symmetry}
To begin with, let us recall a few symmetry properties of the abscissae  (\ref{ci}) and of the polynomials (\ref{cebypol}), which we state without proof (they derive from well known properties of Chebyshev polynomials and abscissae):
\begin{equation}\label{simprop}
c_i = 1-c_{s-i+1}, \quad i=1,\dots,s, \qquad P_j(1-c) = (-1)^jP_j(c), \quad  c\in[0,1], \quad i=0,\dots,s.
\end{equation}
For later use, we also set the vector
\begin{equation}\label{Is1}
\I_s(1) ~:=~ \pmatrix{c} \int_0^1 P_0(x)\dd x\\ \vdots \\ \int_0^1 P_{s-1}(x)\dd x\endpmatrix ~\equiv~
\pmatrix{c} 1\\ 0\\ \frac{-\sqrt{2}}{2^3-1} \\ 0 \\ \frac{-\sqrt{2}}{4^3-1}\\ 0 \\ \frac{-\sqrt{2}}{6^3-1} \\ \vdots\endpmatrix,
\end{equation}
with the last equality following from (\ref{intP1}), such that, according to what seen in the proof of Theorem~\ref{buttabc}:
\begin{equation}\label{bt}
\bfb^\top = \I_s(1)^\top \P_s^\top\Omega.
\end{equation}
We also need to define the following matrices:
\begin{equation}\label{ePD}
P = \pmatrix{ccc} & &1\\ &\udots\\ 1\endpmatrix, \quad D= \pmatrix{cccccc} 1\\ &-1\\ &&1\\&&&-1\\&&&&\ddots\\&&&&&(-1)^{s-1}\endpmatrix \in\RR^{s\times s}.
\end{equation}
We can now state the following result.

\begin{theo} The method  (\ref{tabCCM})-(\ref{bc}) is symmetric.\end{theo}

\proof
Since $P\bfc = \bfc$, see (\ref{simprop}), it is known that the symmetry of the method (\ref{tabCCM})-(\ref{bc}) is granted, provided that \cite{GNI2006}:
\begin{equation}\label{simme}
P\bfb = \bfb, \qquad P\left(\I_s\P_s^\top\Omega\right) P = \bfe\,\bfb^\top-\I_s\P_s^\top\Omega.
\end{equation}
Let us start proving the first equality which, in turn, amounts to show that, by virtue of (\ref{bt}):
$$\bfb^\top  ~=~  \I_s(1)^\top \P_s^\top\Omega P ~=~ \bfb^\top P.$$
In fact, by taking into account (\ref{PIO}), (\ref{simprop}), and (\ref{ePD}), one at first derives that:
\begin{equation}\label{simPO}
\P_s^\top\Omega P = D\P_s^\top\Omega.
\end{equation}
Moreover, by considering that $D\I_s(1)=\I_s(1)$, since  the even entries of $\I_s(1)$ are zero (see (\ref{Is1})), one has:
$$ 
\bfb^\top ~=~  \I_s(1)^\top \P_s^\top\Omega ~=~  \I_s(1)^\top D \P_s^\top\Omega
               ~=~  \I_s(1)^\top \P_s^\top\Omega P ~=~ \bfb^\top P.
$$ 
Consequently, the first part of the statement follows. Further, one obtains, by virtue of (\ref{simPO}):
$$P(\I_s\P_s^\top\Omega) P = P(\I_s D^2 \P_s^\top\Omega) P = (P\I_s D)(D\P_s^\top\Omega P)
=(P\I_s D)\P_s^\top\Omega .$$ 
Therefore, from (\ref{bt}) it follows that the second statement in (\ref{simme}) holds true, provided that
$$P\I_s D = \bfe\,\I_s(1)^\top-\I_s \qquad \Longleftrightarrow\qquad (P\I_sD)_{ij} = \int_{c_i}^1 P_{j-1}(\xi)\dd \xi,\quad i,j=1,\dots,s.$$
From (\ref{PIO}) and (\ref{ePD}) one obtains, by virtue of (\ref{simprop}):
\begin{eqnarray*}
 (P\I_sD)_{ij}  &=& (-1)^{j-1}\int_0^{c_{s-i+1}}P_{j-1}(x)\dd x ~=~\int_0^{1-c_i}(-1)^{j-1}P_{j-1}(x)\dd x\\
 &=& \int_0^{1-c_i}P_{j-1}(1-x)\dd x ~=~ -\int_1^{c_i} P_{j-1}(\xi)\dd\xi ~=~\int_{c_i}^1 P_{j-1}(\xi)\dd\xi.
\end{eqnarray*}
Consequently, the statement is proved.\,\QED\bigskip

\subsection{Stability}\label{Astab} Because of its symmetry,  one concludes that for all $s\ge1$ the absolute region of a CCM$(s)$ method coincides with $\CC^-$, the negative-real complex plane, provided that all the eigenvalues of $X_s$ have positive real part (in fact, $X_s$ is similar to the Butcher matrix, see (\ref{tabCCM})). We have numerically verified that $X_s$ has all the eigenvalues with positive real part for all $s\ge1$ (actually, we stopped at $s=1000$).
We can then conclude that CCM$(s)$ methods are {\em perfectly (or precisely)} $A$-stable for all $s\ge1$.

\subsection{Order of convergence}

Let us now study the order of convergence of a CCM$(s)$ method. For this purpose, we need the following preliminary results.

\begin{lem}\label{symeven} The convergence order of a symmetric method is even.\end{lem}

\proof See \cite[Theorem~3.2]{GNI2006}.\,\QED\bigskip

\begin{lem}\label{Ghj} Assume that a function $G:[0,h]\rightarrow V$, with $V$ a vector space, admits a Taylor expansion at 0. Then, with reference to the orthonormal basis (\ref{cebyw})-(\ref{cebypol}), one has:
$$\psi_j ~:=~\int_0^1 \omega(c) P_j(c)G(ch)\dd c ~=~ O(h^j), \qquad j=0,1,\dots.$$
\end{lem}
\proof In fact,  one has:
\begin{eqnarray*}
\psi_j &=& \int_0^1 \omega(c) P_j(c) G(ch)\dd c ~=~ \int_0^1 \omega(c) P_j(c) \sum_{\ell\ge0}\frac{G^{(\ell)}(0)}{\ell!} (ch)^\ell \dd c \\[1mm]
         &=& \sum_{\ell\ge0} \frac{G^{(\ell)}(0)}{\ell!} h^\ell \int_0^1 \omega(c) P_j(c)c^\ell \dd c ~=~\sum_{\ell\ge j} \frac{G^{(\ell)}(0)}{\ell!} h^\ell \int_0^1 \omega(c) P_j(c)c^\ell \dd c\\[1mm]
         &=&O(h^j).\,\QED
\end{eqnarray*}\smallskip

\begin{lem}\label{hgj} With reference to the approximate Fourier coefficients defined in (\ref{gamis}),  one has:  $$\hat\gamma_j=O(h^j), \qquad j=0,\dots,s-1.$$\end{lem}
\proof
In fact, from the previous Lemma~\ref{Ghj} applied to $G(ch)=f(u(ch))$, we know that (using the notation (\ref{gamiws}))
\begin{equation}\label{bgamj}
\gamma_j(u) ~=~ \int_0^1\omega(c)P_j(c)f(u(ch))\dd c ~=~ O(h^j).
\end{equation}
On the other hand, the quadrature error,
\begin{equation}\label{Dhj}
\Delta_j(h) ~:=~ \gamma_j(u)-\hat\gamma_j ~=~ O(h^{2s-j}).
\end{equation}
Consequently,
$$\hat\gamma_j ~=~ \gamma_j(u) -\Delta_j(h) ~=~ O(h^j), \qquad j=0,\dots,s-1.\QED$$\smallskip

Finally, we need to recall some well-known perturbation results for ODE-IVPs. For this purpose, let us denote
~$y(t,\xi,\eta)$~ the solution of the problem (compare with (\ref{ode1})),
$$\dot y = f(y), \qquad t>\xi, \qquad y(\xi)=\eta\in\RR^m.$$
Then:
$$
\frac{\partial}{\partial t}y(t,\xi,\eta)~=~ f(y(t,\xi,\eta)),\qquad
\frac{\partial}{\partial \xi}y(t,\xi,\eta) ~=~ -\Phi(t,\xi,\eta) f(\eta),
$$
having set 
$$\Phi(t,\xi,\eta) ~\equiv~ \frac{\partial}{\partial \eta}y(t,\xi,\eta),$$
the solution of the variational problem
$$\dot\Phi(t,\xi,\eta) ~=~ f'(y(t,\xi,\eta))\Phi(t,\xi,\eta), \qquad t>\xi, \qquad \Phi(\xi,\xi,\eta)=I_m.$$

We can now state the convergence result.

\begin{theo}\label{ords} With reference to the polynomial approximation (\ref{uch}) defined by the CCM$(s)$ method (\ref{tabCCM})-(\ref{bc}) to the solution of (\ref{ode1}), one has:
$$\|u-y\| = O(h^{s+1}), \qquad y_1-y(h) = O(h^{r+1}),$$
with ~$r=s$~ if ~$s$~ is even, or ~$r=s+1$,~ otherwise.
\end{theo}
\proof By taking into account the above arguments, and using the notation (\ref{bgamj})-(\ref{Dhj}), one has:
\begin{eqnarray*}
u(ch)-y(ch)&=& y(ch,ch,u(ch))-y(ch,0,u(0)) ~=~ \int_0^{ch} \frac{\dd}{\dd t} y(ch,t,u(t))\,\dd t\\[1mm]
&=& \int_0^{ch}\left.\frac{\partial}{\partial \xi}y(ch,\xi,u(t))\right|_{\xi=t}+\left.\frac{\partial}{\partial \eta}y(ch,t,\eta)\right|_{\eta=u(t)}\dot u(t)\,\dd t \\
\end{eqnarray*}\begin{eqnarray*}
&=& h\int_0^{c}\left.\frac{\partial}{\partial \xi}y(ch,\xi, u(\tau h))\right|_{\xi=\tau h}+\left.\frac{\partial}{\partial \eta}y(ch,\tau h,\eta)\right|_{\eta=u(\tau h)}\dot u(\tau h)\,\dd \tau \\[1mm]
&=&-h\int_0^c \Phi(ch,\tau h,u(\tau h))\left[ f(u(\tau h))-\sum_{j=0}^{s-1}P_j(\tau)\hat\gamma_j\right]\dd\tau \\[1mm]
&=&-h\int_0^c \Phi(ch,\tau h,u(\tau h))\left[ \sum_{j\ge0} P_j(\tau)\gamma_j(u)-\sum_{j=0}^{s-1}P_j(\tau)\Big(\gamma_j(u)-\Delta_j(h)\Big)\right]\dd\tau \\[1mm]
&=&-h\int_0^c \Phi(ch,\tau h,u(\tau h))\left[ \sum_{j\ge s} P_j(\tau)\gamma_j(u)+\sum_{j=0}^{s-1}P_j(\tau)\Delta_j(h)\right]\dd\tau \\[1mm]
&=&-h\int_0^c \Phi(ch,\tau h,u(\tau h))\left[O(h^s)+O(h^{s+1})\right]\dd\tau ~=~ O(h^{s+1}).
\end{eqnarray*}
Consequently, the first part of the statement holds. 
The last part of the statement then follows from Lemma~\ref{symeven}, by taking $c=1$.\,\QED\bigskip

\subsection{Analysis as spectral method}

As was pointed out in the introduction, one main reason for introducing CCMs is their use as spectral methods in time. This strategy allows for using large (sometimes huge) timesteps $h$. In this context, the classical notion of order of convergence, based on the fact that $h\rightarrow0$,  does not apply anymore, so that a different analysis is needed. We shall here follow similar steps as those described in \cite{ABI2020} for spectral HBVMs. 

To begin with, we recall that, for the modified Chebysehv basis (\ref{cebyw})-(\ref{cebypol}) the following properties hold true:\footnote{The second property has already been used in the proof of Lemma~\ref{Ghj}.}
\begin{equation}\label{prop}
\forall j\ge0:\quad \|P_j\|\le\sqrt{2}, \qquad \int_0^1 \omega(c)P_j(c)c^\ell\dd c = 0, \quad  \ell=0,\dots,j-1.
\end{equation}
Moreover, we recall that, with reference to the notation (\ref{gamiws}):\footnote{Also this property has already been used, in the proof of Theorem~\ref{ords}.}
\begin{equation}\label{fych}
f(y(ch)) = \sum_{j\ge0} P_j(c)\gamma_j(y), \qquad c\in[0,1].
\end{equation}
We also state the following preliminary result.

\begin{lem}\label{Qj}
With reference to (\ref{cebyw}) and to the polynomials (\ref{cebypol}), let us define, for a given $\rho>1$:
$$Q_j(\xi) ~=~ \int_0^1 \omega(c)\frac{P_j(c)}{\xi-c}\dd c, \qquad \xi\in\C_\rho(0),$$
having set ~$\C_\rho(0)=\{z\in\CC:|z|=\rho\}$. Then
\begin{equation}\label{normro}
\|Q_j\|_\rho ~:=~ \max_{\xi\in\C_\rho(0)}|Q_j(\xi)| ~\le~ \sqrt{2}\frac{\rho^{-j}}{(\rho-1)}.
\end{equation}
\end{lem}
\proof By taking into account (\ref{prop}), for $|\xi|=\rho>1$ one has:
\begin{eqnarray*}
Q_j(\xi) &=&\int_0^1 \omega(c)\frac{P_j(c)}{\xi-c}\dd c~=~\xi^{-1} \int_0^1 \omega(c)\frac{P_j(c)}{1-\xi^{-1}c}\dd c\\[1mm]
&=& \xi^{-1} \int_0^1 \omega(c)P_j(c)\sum_{\ell\ge0}\xi^{-\ell}c^\ell\dd c ~=~
\xi^{-1}\sum_{\ell\ge j}\xi^{-\ell} \int_0^1 \omega(c)P_j(c)c^\ell\dd c\\[1mm]
&=& \xi^{-j-1}\sum_{\ell\ge 0}\xi^{-\ell} \int_0^1 \omega(c)P_j(c)c^{\ell+j}\dd c.
\end{eqnarray*}
Evaluating the norms, one has, again by virtue of (\ref{prop}):
$$\left|\int_0^1 \omega(c)P_j(c)c^{\ell+j}\dd c\right| ~\le~ \|P_j\| \int_0^1 \omega(c)\dd c ~\le~ \sqrt{2}.$$
Consequently, one obtains:
$$|Q_j(\xi)|~\le~ \sqrt{2}\rho^{-j-1}\sum_{\ell\ge0} \rho^{-\ell} ~=~ \sqrt{2}\rho^{-j-1}\frac{1}{1-\rho^{-1}} ~=~ \sqrt{2}\frac{\rho^{-j}}{\rho-1}.\,\QED$$\smallskip

Further, by setting hereafter, for $r>0$,
$$\B_r(0) = \{z\in\CC:|z|\le r\},$$
we have the following straightforward result.

\begin{lem}\label{ghlem}
Let $g(z)$ be analytical in $\B_{r^*}(0)$, for a given $r^*>0$. Then, for all $h\in(0,h^*]$, with $h^*<r^*$, 
\begin{equation}\label{ghz}
g_h(\xi) ~:=~ g(\xi h)
\end{equation}
is analytical in $\B_\rho(0)$, with: 
\begin{equation}\label{rostar}
\rho ~=~ \rho(h) ~:=~ \frac{r^*}h ~\ge~  \frac{r^*}{h^*} ~=:~ \rho^*~>~1.
\end{equation}
\end{lem}

We can now state the following result.

\begin{theo}\label{grhoj}
With reference to (\ref{ode1}), let assume that the function
\begin{equation}\label{gdiz}
g(z)~:=~f(y(z))
\end{equation}
and the timestep $h$ satisfy the hypotheses of Lemma~\ref{ghlem}. Then there exists $\kappa=\kappa(h^*)$ and $\rho>1$, $\rho\sim h^{-1}$, such that
\begin{equation}\label{gjkro}
|\gamma_j(y)|~\le~ \kappa \rho^{-j}, \qquad j\ge0.
\end{equation}
\end{theo}
\proof By taking into account (\ref{gamiws}), (\ref{fych}), and (\ref{ghz})-(\ref{rostar}), by virtue of Lemma~\ref{Qj} one has:
\begin{eqnarray*}
\gamma_j(y) &=& \int_0^1 \omega(c)P_j(c)f(y(ch))\dd c ~\equiv~\int_0^1 \omega(c)P_j(c)g_h(c)\dd c\\[1mm]
&=&\int_0^1 \omega(c)P_j(c) \left[\frac{1}{2\pi i}\int_{\C_\rho(0)} \frac{g_h(\xi)}{\xi-c}\dd\xi\right] \dd c 
~=~ \frac{1}{2\pi i}\int_{\C_\rho(0)} g_h(\xi) \left[\int_0^1\omega(c)\frac{P_j(c)}{\xi-c}\dd c \right]\dd \xi \\[1mm]
&\equiv& \frac{1}{2\pi i}\int_{\C_\rho(0)} g_h(\xi) Q_j(\xi)\dd \xi.
\end{eqnarray*}
Consequently (see (\ref{normro})),
$$|\gamma_j(y)|~\le~ \frac{\rho}{2\pi} \|g_h\|_\rho\|Q_j\|_\rho.$$
From Lemma~\ref{Qj} we know that 
$$\|Q_j\|_\rho ~\le~ \sqrt{2} \frac{\rho^{-j}}{\rho-1}.$$ 
Further,
$$\|g_h\|_\rho ~=~ \max_{|z|=\rho}|g_h(z)| ~\le~  \max_{|z|\le \rho^*}|g_h(z)| ~\equiv~ 
 \max_{|z|\le r^*}|g(z)| ~=:~ \|g\|.$$
 Consequently, one eventually derives:
 $$|\gamma_j(y)|~\le~ \frac{\rho\sqrt{2}\|g\|}{2\pi(\rho-1)} \rho^{-j} ~\le~ 
 \overbrace{\frac{\rho^*\sqrt{2}\|g\|}{2\pi(\rho^*-1)}}^{=:\,\kappa} \rho^{-j},
 $$
 from which the statement follows by taking into account that (see (\ref{rostar})), $\rho^*=r^*/h^*$.\,\QED\bigskip

Finally, we have the following result.

\begin{theo}\label{valeancora}
Let us consider the polynomial approximation (\ref{sig}) and (\ref{cebyw})-(\ref{cebypol}), with the Fourier coefficients given by (\ref{gamiws}). Then, Theorem~\ref{spectral} continues formally to hold.
\end{theo}

\medskip
\proof Following similar steps as in the proof of Theorem~\ref{ords}, one has:
\begin{eqnarray*}
\sigma(ch)-y(ch)&=& y(ch,ch,\sigma(ch))-y(ch,0,\sigma(0)) ~=~ \int_0^{ch} \frac{\dd}{\dd t} y(ch,t,\sigma(t))\,\dd t\\[1mm]
&=& \int_0^{ch}\left.\frac{\partial}{\partial \xi}y(ch,\xi,\sigma(t))\right|_{\xi=t}+\left.\frac{\partial}{\partial \eta}y(ch,t,\eta)\right|_{\eta=\sigma(t)}\dot \sigma(t)\,\dd t \\[1mm]
\hspace{2.2cm}
&=& h\int_0^{c}\left.\frac{\partial}{\partial \xi}y(ch,\xi, \sigma(\tau h))\right|_{\xi=\tau h}+\left.\frac{\partial}{\partial \eta}y(ch,\tau h,\eta)\right|_{\eta=\sigma(\tau h)}\dot \sigma(\tau h)\,\dd \tau \\[1mm]
&=&-h\int_0^c \Phi(ch,\tau h,\sigma(\tau h))\left[ f(\sigma(\tau h))-\sum_{j=0}^{s-1}P_j(\tau)\gamma_j(\sigma)\right]\dd\tau \\[1mm]
\end{eqnarray*}\begin{eqnarray*}
&=&-h\int_0^c \Phi(ch,\tau h,\sigma(\tau h))\left[ \sum_{j\ge0} P_j(\tau)\gamma_j(\sigma)-\sum_{j=0}^{s-1}P_j(\tau)\gamma_j(\sigma)\right]\dd\tau \\[1mm]
&=&-h\int_0^c \Phi(ch,\tau h,\sigma(\tau h))\sum_{j\ge s} P_j(\tau)\gamma_j(\sigma)\dd\tau.
\end{eqnarray*}
Consequently, from (\ref{prop}) and Theorem~\ref{grhoj}, also considering (\ref{rostar}), one derives:
\begin{eqnarray*}
\|\sigma-y\| &\le& h\sqrt{2}\kappa \max_{c,\tau\in[0,1]}\|\Phi(ch,\tau h,\sigma(\tau h))\|  \sum_{j\ge s}\rho^{-j} \\
                   &\le& h\underbrace{\sqrt{2}\kappa \max_{c,\tau\in[0,1]}\|\Phi(ch^*,\tau h^*,\sigma(\tau h^*))\| \frac{1}{1-(\rho^*)^{-1}}}_{=:\,M}\rho^{-s} ~\equiv~ hM\rho^{-s}.\,\QED
\end{eqnarray*}\medskip

\begin{rem}\label{uedsig}
When a CCM$(s)$ method is used as a spectral method in time, one obviously assume that the value of $s$ is large enough so that the quadrature errors falls below the round-off error level. In other words, the polynomials $\sigma$ and $u$ are undistinguishable, within the round-off error level.
\end{rem}

\section{Numerical tests}\label{numtest}

We here report a few numerical tests to assess the theoretical findings. They have been carried out on a 3 GHz Intel Xeon W10 core computer with 64GB of memory, running Matlab$^\copyright$ ~2020b.  We consider the Kepler problem,

\begin{equation}\label{kepl}
\dot q_i = p_i, \qquad \dot p_i = \frac{-q_i}{(q_1^2+q_2^2)^{\frac{3}2}}, \qquad i=1,2,
\end{equation}
that, when considering the trajectory starting at
\begin{equation}\label{kepl0}
q_1(0) = 0.4, \qquad q_2(0) = p_1(0) = 0, \qquad p_2(0)=2,
\end{equation}
has a periodic solution of period $2\pi$. In Table~\ref{tab1} we list the numerical results obtained by solving this problem on one period by means of the CCM$(s)$ method, $s=1,\dots,4$, using timestep $h=2\pi/n$, namely, the errorr {\em err} after one period and the corresponding rate of convergence. As one may see, the listed results agree with what stated in Theorem~\ref{ords}, in particular the fact that the order of the methods is even, due to their symmetry.

\begin{table}[t]
\caption{Numerical results when solving problem (\ref{kepl})-(\ref{kepl0}) on the interval $[0,2\pi]$ with timestep $h=2\pi/n$.}\label{tab1}
\centering
\begin{tabular}{|r|cc|cc|cc|cc|}
\hline
   &\multicolumn{2}{c|}{$s=1$}&\multicolumn{2}{c|}{$s=2$}&\multicolumn{2}{c|}{$s=3$}&\multicolumn{2}{c|}{$s=4$}\\
   \hline
 $n$ & \em err & rate &  \em err & rate &\em err & rate & \em err & rate \\
\hline
   50 & 2.98e+0 &   --- & 2.24e+0 &   --- & 7.36e-03 &   --- & 7.33e-03 &   --- \\ 
  100 & 1.66e+0 &   0.8 & 9.45e-01 &   1.2 & 6.15e-04 &   3.6 & 4.46e-04 &   4.0 \\ 
  200 & 5.23e-01 &   1.7 & 2.53e-01 &   1.9 & 4.03e-05 &   3.9 & 2.78e-05 &   4.0 \\ 
  400 & 1.34e-01 &   2.0 & 6.34e-02 &   2.0 & 2.55e-06 &   4.0 & 1.73e-06 &   4.0 \\ 
  800 & 3.35e-02 &   2.0 & 1.58e-02 &   2.0 & 1.60e-07 &   4.0 & 1.08e-07 &   4.0 \\ 
 1600 & 8.38e-03 &   2.0 & 3.96e-03 &   2.0 & 1.00e-08 &   4.0 & 6.77e-09 &   4.0 \\ 
\hline
\end{tabular}

\bigskip
\caption{Numerical results when solving problem (\ref{kepl})-(\ref{kepl0}) on the interval $[0,20\pi]$ by using the CCM(50) method with timestep $h=2\pi/n$.}\label{tab2}
\centering
\begin{tabular}{|c|c|c|c|c|c|}
\hline
  $n$ & 3 & 6 & 9 & 12 & 15 \\
   \hline
 period & \em err & \em err & \em err & \em err & \em err \\
\hline
  1 & 5.04e-12 & 9.14e-14 & 7.37e-14 & 1.27e-13 & 7.44e-14 \\ 
  2 & 9.72e-12 & 7.96e-14 & 1.24e-13 & 3.05e-13 & 5.72e-14 \\ 
  3 & 1.34e-11 & 3.54e-13 & 2.81e-13 & 7.23e-13 & 5.49e-14 \\ 
  4 & 1.90e-11 & 3.00e-13 & 7.80e-13 & 1.27e-12 & 1.49e-13 \\ 
  5 & 2.55e-11 & 4.69e-13 & 1.34e-12 & 2.30e-12 & 2.77e-13 \\ 
  6 & 3.04e-11 & 5.68e-13 & 1.75e-12 & 3.25e-12 & 3.59e-13 \\ 
  7 & 3.44e-11 & 2.13e-13 & 1.73e-12 & 4.23e-12 & 2.37e-13 \\ 
  8 & 3.94e-11 & 2.72e-13 & 1.45e-12 & 5.19e-12 & 2.04e-13 \\ 
  9 & 4.38e-11 & 6.93e-13 & 1.18e-12 & 6.15e-12 & 3.29e-13 \\ 
 10 & 4.77e-11 & 1.54e-12 & 8.38e-13 & 7.01e-12 & 5.00e-13 \\ 
\hline
\end{tabular}

\end{table}

\begin{figure}[t]
\centerline{
\includegraphics[width=12cm]{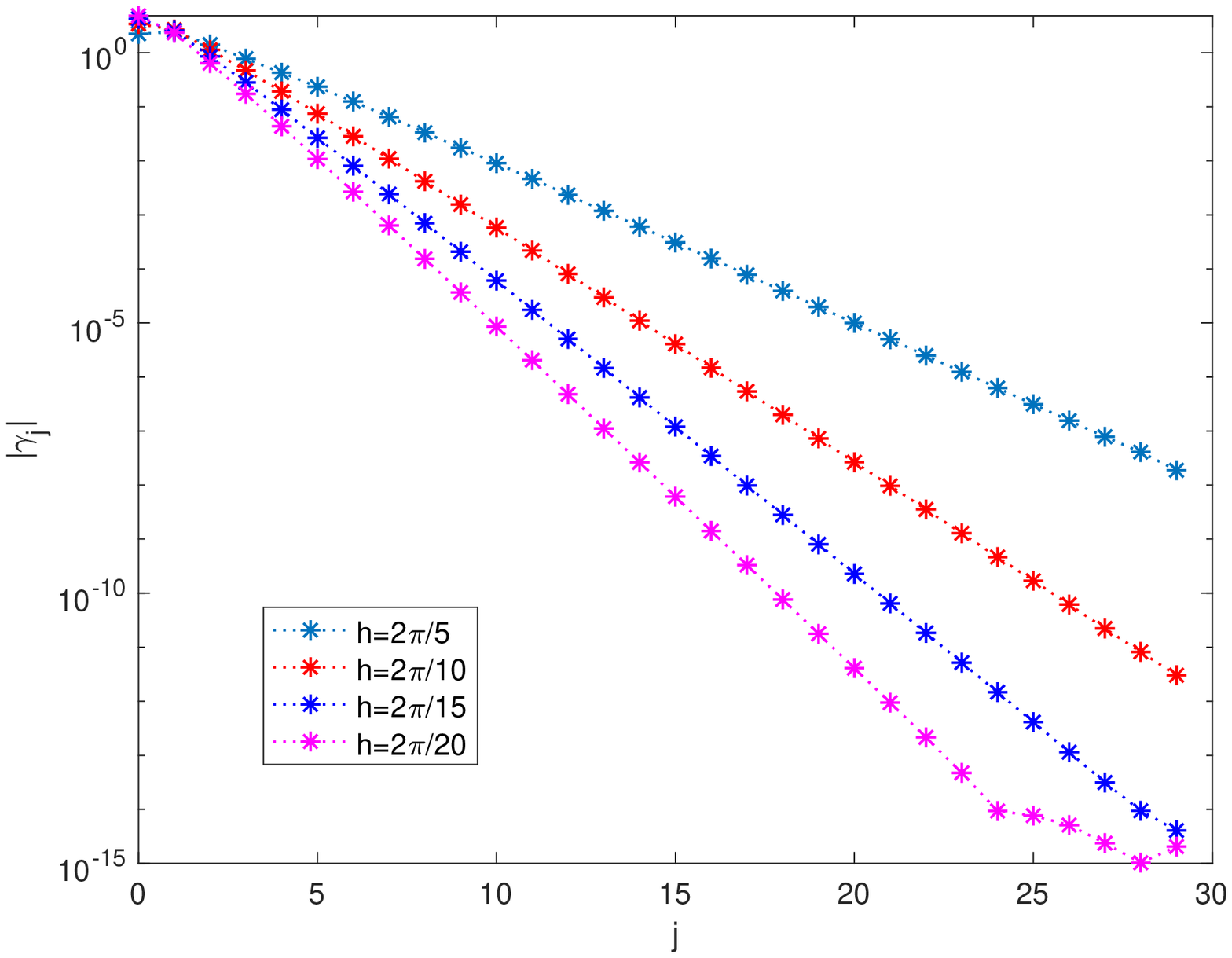}}
\caption{Plot of $|\gamma_j|$, $j=0,\dots,29$, for the CCM$(30)$ method, for the first integration step of problem (\ref{kepl})-(\ref{kepl0}) by using the given timestep $h$.}\label{fig1}
\end{figure}

Next, we consider the use of the methods as a spectral methods in time. In Figure~\ref{fig1} we plot the values of $|\gamma_j|$, $j=0,\dots,29$, for the CCM(30) method used with timestep $h=2\pi/n$, $n=5,10,15,20$, for the first integration step. As one may see, the Fourier coefficients decrease exponentially, with $j$, and the basis of the exponential decreases with $h$ (i.e., they decrease as $\rho^j$, with $\rho\sim h^{-1}$). Moreover, when using $h=\pi/10$ the last Fourier coefficients reach the round-off error level, thus stagnating.

Now, let us solve problem (\ref{kepl})-(\ref{kepl0}) on ten periods, by using the CCM(50) method with timestep $h=2\pi/n$, $n=3,6,9,12,15$. The errors ({\em err}), at the end of each period, are listed in Table~\ref{tab2}. As one my see, the errors are uniformly small, as is expected from a spectral method.

\begin{figure}[t]
\centerline{
\includegraphics[width=7.5cm]{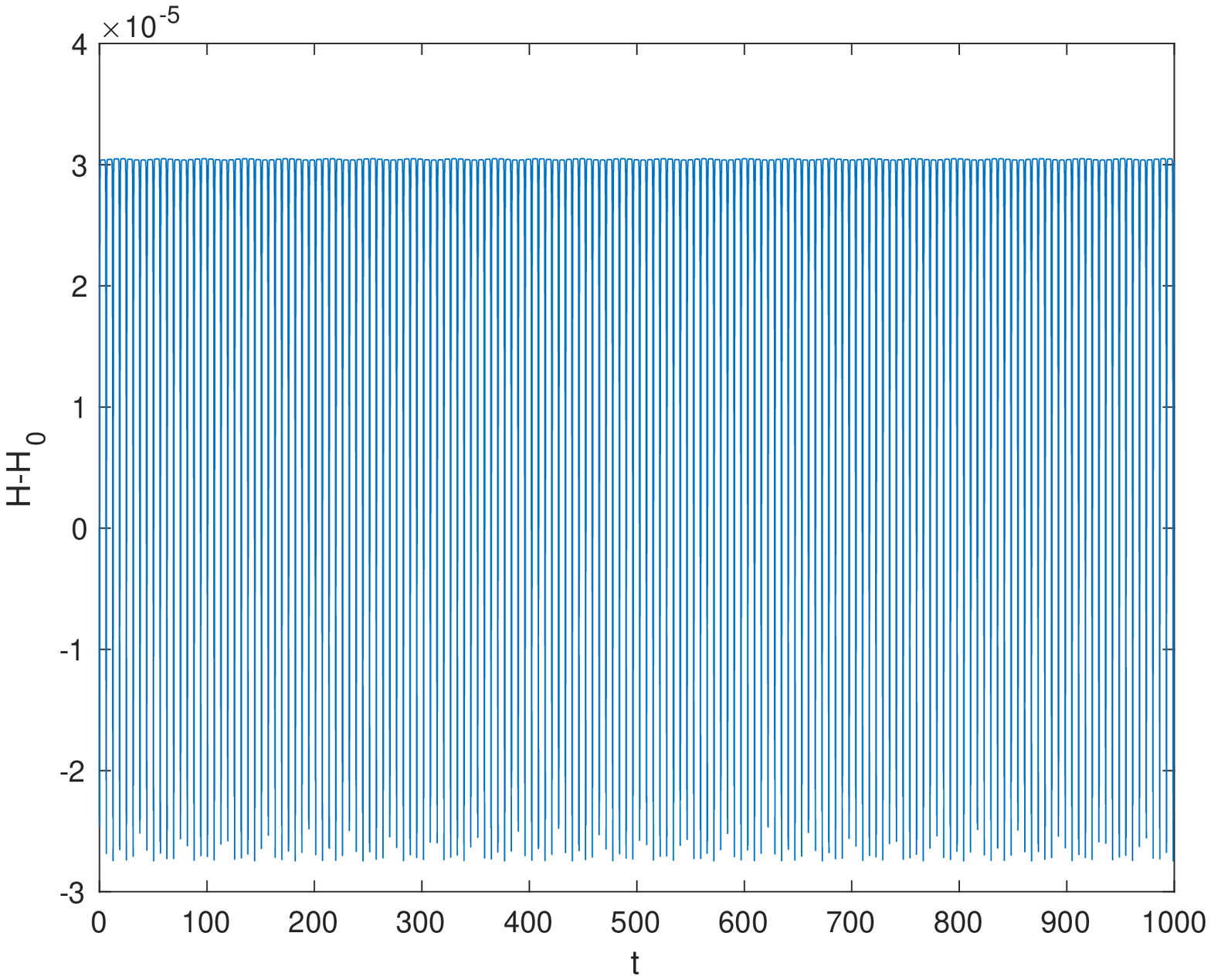}\quad
\includegraphics[width=7.5cm]{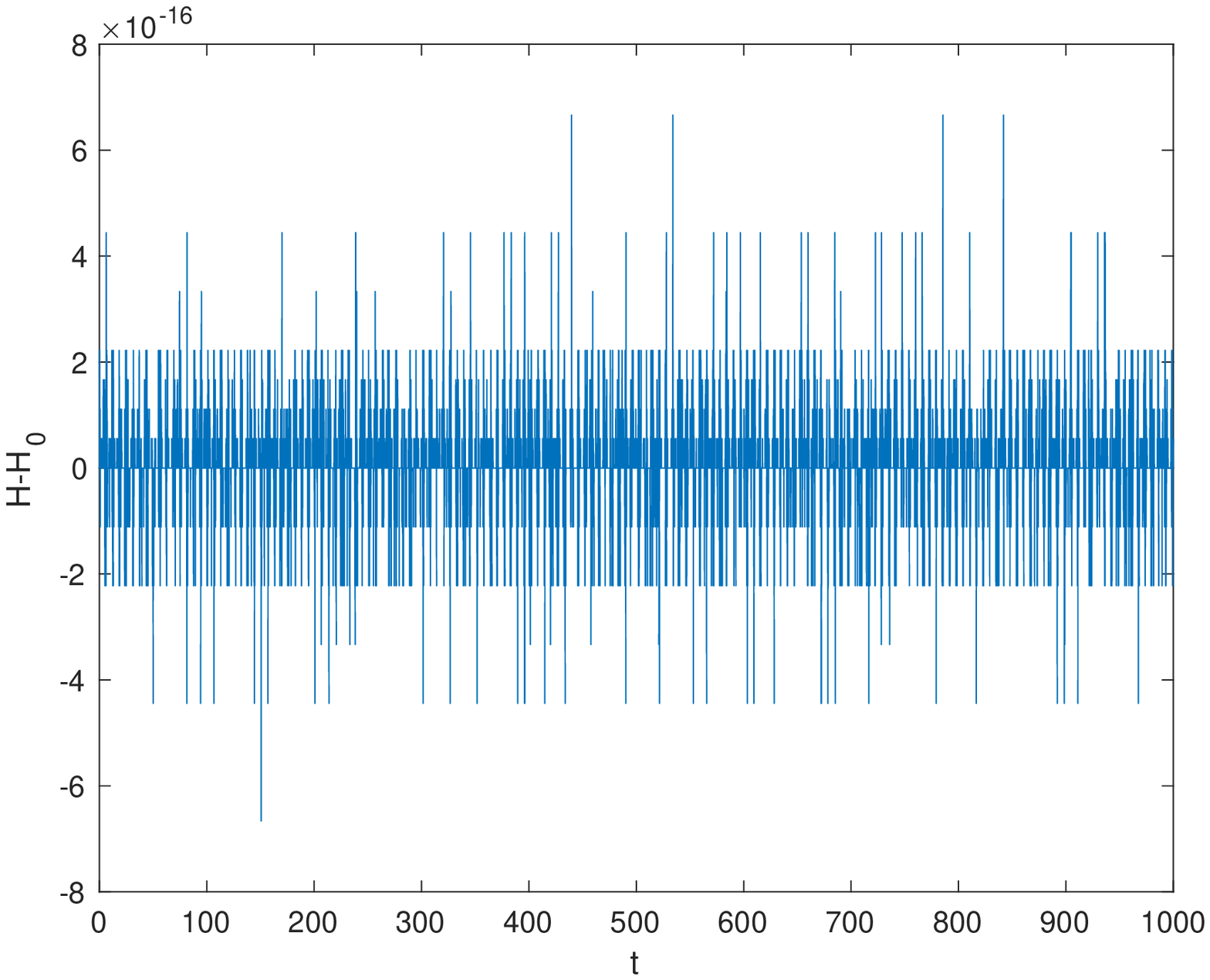}}
\caption{Plot of the Hamiltonian error when solving the problem (\ref{kepl})-(\ref{kepl0}) with timestep $h=0.1$ one the interval $[0,10^3]$ by using the CCM(3) method (left plot) and the CCM(30) method (right plot).}\label{fig2}
\end{figure}

At last, since problem (\ref{kepl})-(\ref{kepl0}) is Hamiltonian with Hamiltonian
$$H(q,p) = \frac{1}2 \|p\|_2^2 - \frac{1}{\|q\|_2},$$
let us consider the solution of problem (\ref{kepl})-(\ref{kepl0}) on the interval $[0,10^3]$ with timestep $h=0.1$, by using the CCM(3) and CCM(30) methods. The former method is only symmetric, whereas the second one is spectral, for the considered timestep. On the left of Figure~\ref{fig2} is the plot of the Hamiltonian error for the CCM(3) method: though not energy conserving, the Hamiltonian error is bounded, since the method is symmetric  and the problem reversible. On the right of the same figure, is the plot of the Hamiltonian error for the CCM(30) method: in this case, the spectral accuracy reflects on the practical conservation of the Hamiltonian (as well as of the additional constant of motions, such as the angular momentum and the Lenz vector, not displayed here). It must be emphasized that the execution times for running the two methods (implemented in the same code) is almost the same: 2.9 sec vs. 3.7 sec. Considering that a similar behavior is observed when solving other classes of problems (see also  \cite{ABI2020,BIMR2019} for a related analysis on HBVMs), it is then clear that the use of CCMs as spectral methods can be very effective.

\section{Conclusions}\label{fine}

In this paper we have re-derived, in a novel framework, the Chebyshev collocation methods of Costabile and Napoli \cite{CoNa2001}, thus providing a more comprehensive analysis of the methods, w.r.t. \cite{CoNa2001,CoNa2004}, in particular when used as spectral methods in time. The methods, derived by considering an expansion of the vector field along the Chebyshev orthonormal polynomial basis, turn out to be symmetric, perfectly $A$-stable, can have arbitrarily high order, and the entries of their Butcher tableau can be written in closed form. Their use  as spectral methods in time has been also studied. A few numerical tests confirm the theoretical achievements. As a further direction of investigation, the efficient implementation of the methods will be considered, as well as their application to different kinds of differential problems.

\subsection*{Declarations}

The authors declare no conflict of interests.

\subsection*{Acknoledgements} 

The authors wish to thank the {\tt mrSIR} project \cite{mrSIR} for the financial support.


\begin{thebibliography}{99}

\bibitem{ABI2015} P.\,Amodio, L.\,Brugnano, F.\,Iavernaro. Energy-conserving methods for Hamiltonian Boundary Value Problems and applications in astrodynamics. {\em  Adv. Comput. Math.} {\bf 41} (2015) 881--905. \url{https://doi.org/10.1007/s10444-014-9390-z}

\bibitem{ABI2019_0} P.\,Amodio,  L.\,Brugnano, F.\,Iavernaro. A note on the continuous-stage Runge-Kutta-(Nystr\"om) formulation of Hamiltonian Boundary Value Methods (HBVMs).  {\em Appl. Math. Comput.} {\bf 363} (2019) 124634.  \url{https://doi.org/10.1016/j.amc.2019.124634}

\bibitem{ABI2022_0} P.\,Amodio,  L.\,Brugnano, F.\,Iavernaro. Continuous-Stage Runge-Kutta approximation to Differential Problems. {\em Axioms} {\bf 11} (2022) 192. \url{https://doi.org/10.3390/axioms11050192} 

\bibitem{ABI2019} P.\,Amodio, L.\,Brugnano, F.\,Iavernaro. Spectrally accurate solutions of nonlinear fractional initial value problems. {\em AIP Conf. Proc.} {\bf 2116} (2019) 140005. \url{https://doi.org/10.1063/1.5114132}

\bibitem{ABI2020} P.\,Amodio, L.\,Brugnano, F.\,Iavernaro. Analysis of Spectral Hamiltonian Boundary Value Methods (SHBVMs) for the numerical solution of ODE problems.  {\em Numer. Algorithms}   {\bf 83} (2020) 1489--1508. \url{https://doi.org/10.1007/s11075-019-00733-7}

\bibitem{ABI2022} P.\,Amodio, L.\,Brugnano, F.\,Iavernaro.  Arbitrarily high-order energy-conserving methods for Poisson problems. {\em Numer. Algorithms} (2022). \url{https://doi.org/10.1007/s11075-022-01285-z}

\bibitem{ABIM2020}  P.\,Amodio, L.\,Brugnano, F.\,Iavernaro, C.\,Magherini. Spectral solution of ODE-IVPs using SHBVMs. {\em AIP Conf. Proc.} {\bf 2293} (2020) 100002. \url{https://doi.org/10.1063/5.0026404}

\bibitem{BBTZ2020} L.\,Barletti, L.\,Brugnano, Y.\,Tang, B.\,Zhu. Spectrally accurate space-time solution of Manakov systems.  {\em J.~Comput. Appl. Math.} {\bf  377} (2020) 112918.
\url{https://doi.org/10.1016/j.cam.2020.112918}

\bibitem{BlCa2016} S.\,Blanes, F.\,Casas. {\em A Concise Introduction to Geometric Numerical Integration}. CRC Press, Boca Raton, FL, USA, 2016.

\bibitem{BCMR2012} L.\,Brugnano, M.\,Calvo, I.J.\,Montijano, L.\,R\`andez.  Energy preserving methods for Poisson systems. {\em J.~Comput. Appl. Math.} {\bf  236} (2012) 3890--3904. \url{http://doi.org/10.1016/j.cam.2012.02.033}

\bibitem{BFCI2014} L.\,Brugnano, G.\,Frasca-Caccia, F.\,Iavernaro. Efficient  implementation of Gauss collocation and Hamiltonian Boundary Value Methods. {\em Numer. Algorithms} {\bf 65} (2014) 633--650. \url{https://doi.org/10.1007/s11075-014-9825-0}

\bibitem{BGIW2018} L.\,Brugnano, G.\,Gurioli, F.\,Iavernaro, E.\,Weinm\"uller. Line integral solution of Hamiltonian systems with holonomic constraints. {\em Appl. Numer. Math.} {\bf 127} (2018) 56--77. \url{ https://doi.org/10.1016/j.apnum.2017.12.014}

\bibitem{BGZ2019}  L.\,Brugnano, G.\,Gurioli, C.\,Zhang. Spectrally accurate energy-preserving methods for the numerical solution of the ``Good'' Boussinesq equation. {\em  Numer. Methods Partial Differential Equations}  {\bf 35} (2019) 1343--1362. \url{https://doi.org/10.1002/num.22353}

\bibitem{BI2012} L.\,Brugnano, F.\,Iavernaro.  Line integral methods which preserve all invariants of conservative problems. {\em J.~Comput. Appl. Math.} {\bf  236} (2012)  3905–3919. \url{http://doi.org/10.1016/j.cam.2012.03.026}

\bibitem{LIMBook2016} L.\,Brugnano, F.\,Iavernaro. {\em Line Integral Methods for Conservative Problems}.  Chapman et Hall/CRC, Boca~Raton, FL, USA, 2016.

\bibitem{BI2018}  L.\,Brugnano, F.\,Iavernaro. Line Integral Solution of Differential Problems. {\em Axioms} {\bf 7}(2) (2018) 36. \url{https://doi.org/10.3390/axioms7020036}
 
\bibitem{BIMR2019} L.\,Brugnano, F.\,Iavernaro, J.I.\,Montijano, L.\,R\'andez. Spectrally accurate space-time solution of Hamiltonian PDEs. {\em Numer. Algorithms}   {\bf 81} (2019) 1183--1202.
\url{https://doi.org/10.1007/s11075-018-0586-z }

\bibitem{BIS2010} L.\,Brugnano, F.\,Iavernaro, T.\,Susca. Numerical comparisons between Gauss-Legendre methods and Hamiltonian BVMs defined over Gauss points. {\em Monogr. Real Acad. Cienc. Zaragoza} {\bf 33} (2010) 95--112.

\bibitem{BIT2009}  L.\,Brugnano, F.\,Iavernaro, D.\,Trigiante. Hamiltonian BVMs (HBVMs): a family of  ``drift-free'' methods for integrating polynomial Hamiltonian systems. {\em AIP Conf. Proc.} {\bf 1168} (2009) 715--718.  \url{http://doi.org/10.1063/1.3241566}
 
\bibitem{BIT2010} L.\,Brugnano, F.\,Iavernaro, D.\,Trigiante.  Hamiltonian Boundary Value Methods (Energy Preserving Discrete Line Integral Methods). {\em  JNAIAM J. Numer. Anal. Ind. Appl. Math.} {\bf  5} (2010) 17--37.  

\bibitem{BIT2011} L.\,Brugnano, F.\,Iavernaro, D.\,Trigiante. A note on the efficient implementation of Hamiltonian BVMs. {\em J.~Comput. Appl. Math.} {\bf 236} (2011) 375--383. \url{https://doi.org/10.1016/j.cam.2011.07.022}

\bibitem{BIT2012} L.\,Brugnano, F.\,Iavernaro, D.\,Trigiante. The lack of continuity and the role of infinite and infinitesimal in numerical methods for ODEs: The case of symplecticity.  {\em Appl. Math. Comput.} {\bf 218} (2012) 8053--8063. \url{https://doi.org/10.1016/j.amc.2011.03.022}

\bibitem{BIT2012_1} L.\,Brugnano, F.\,Iavernaro, D.\,Trigiante.  A simple framework for the derivation and analysis of effective one-step methods for ODEs. {\em Appl. Math. Comput.} {\bf 218} (2012) 8475--8485. \url{https://doi.org/10.1016/j.amc.2012.01.074}

\bibitem{BIT2012_2} L.\,Brugnano, F.\,Iavernaro, D.\,Trigiante. A two-step, fourth-order method with energy preserving properties. {\em  Comput. Phys. Commun.}  {\bf 183} (2012) 1860--1868. \url{https://doi.org/10.1016/j.cpc.2012.04.002}

\bibitem{BIT2015} L.\,Brugnano, F.\,Iavernaro, D.\,Trigiante. Analisys of Hamiltonian Boundary Value Methods (HBVMs): A class of energy-preserving Runge-Kutta methods for the numerical solution of polynomial Hamiltonian systems.  {\em Commun. Nonlinear Sci. Numer. Simul.} {\bf 20} (2015) 650--667. \url{https://doi.org/10.1016/j.cnsns.2014.05.030}

\bibitem{BIZ2020}  L.\,Brugnano, F.\,Iavernaro, R.\,Zhang.  Arbitrarily high-order energy-preserving methods for simulating the gyrocenter dynamics of charged particles. {\em J. Comput. Appl. Math.} {\bf 380} (2020) 112994. \url{https://doi.org/10.1016/j.cam.2020.112994}

\bibitem{BMR2019_0} L.\,Brugnano, J.I.\,Montijano, L.\,R\'andez. High-order energy-conserving Line Integral Methods for charged particle dynamics. {\em  J. Comput. Phys.} {\bf 396} (2019) 209--227. \url{https://doi.org/10.1016/j.jcp.2019.06.068}

\bibitem{BMR2019} L.\,Brugnano, J.I.\,Montijano, L.\,R\'andez. On the effectiveness of spectral methods for the numerical solution of multi-frequency highly-oscillatory Hamiltonian problems.  {\em Numer. Algorithms}   {\bf 81} (2019) 345--376. \url{http://doi.org/10.1007/s11075-018-0552-9}

\bibitem{BS2014} L.\,Brugnano, Y.\,Sun. Multiple invariants conserving Runge-Kutta type methods for Hamiltonian problems. {\em Numer. Algorithms}  {\bf 65} (2014) 611--632. \url{https://doi.org/10.1007/s11075-013-9769-9}

\bibitem{But21} J.C.\,Butcher. {\em B-series. Algebraic analysis of numerical methods.} Springer Nature, Switzerland, 2021.

\bibitem{CMcLMcLOQW2009} E.\,Celledoni, R.I.\,McLachlan, D.I.\,McLaren, B.\,Owren, G.R.W.\,Quispel, W.M.\,Wright. Energy-preserving Runge-Kutta methods. {\em M2AN Math. Model. Numer. Anal.} {\bf 43} (2009) 645--649. \url{https://doi.org/10.1051/m2an/2009020}

\bibitem{CFM2006} P.\,Chartier, E.\,Faou, A.\,Murua, An algebraic approach to invariant preserving integrators: the case of quadratic and Hamiltonian invariants. {\em Numer. Math.} {\bf 103} (2006) 575--590. \url{https://doi.org/10.1007/s00211-006-0003-8}

\bibitem{CoNa2001} F.\,Costabile, A.\,Napoli. A method for global approximation of the initial value problem. {\em Numer. Algorithms} {\bf 27 } (2001) 119--130.

\bibitem{CoNa2004} F.\,Costabile, A.\,Napoli. Stability of Chebyshev collocation methods. {\em Comput. Math. Appl.} {\em 47} (2004) 659--666.

\bibitem{CoNa2011} F.\,Costabile, A.\,Napoli. A class of collocation methods for numerical integration of initial value problems. {\em Comput. Math. Appl.}  {\bf 62} (2011) 3221--3235.

\bibitem{H2010} E.\,Hairer. Energy preserving variant of collocation methods. {\em  JNAIAM J. Numer. Anal. Ind. Appl. Math.} {\bf 5} (2010) 73--84.

\bibitem{GNI2006} E.\,Hairer, Ch.\,Lubich, G.\,Wanner. {\em Geometric Numerical Integration, 2nd ed}. Springer, Berlin, Germany, 2006. 

\bibitem{HW1996} E.\,Hairer, G.\,Wanner. {\em Solving ordiary differential equations II: Stiff and differential-algebraic problems, second ed.} Springer-Verlag, Berlin, 1996.

\bibitem{IaPa2007} F.\,Iavernaro, B.\,Pace. $s$-Stage trapezoidal methods for the conservation of Hamiltonian functions of polynomial type. {\em AIP Conf. Proc.}  {\bf 936} (2007) 603--606.
\url{https://doi.org/10.1063/1.2790219}

\bibitem{IaPa2008} F.\,Iavernaro, B.\,Pace. Conservative block-Boundary Value Methods for the solution of polynomial Hamiltonian systems. {\em AIP Conf. Proc.} {\bf 1048} (2008) 888--891.
\url{https://doi.org/10.1063/1.2991075}

\bibitem{IaTr2009} F.\,Iavernaro, D.\,Trigiante. High-order Symmetric Schemes for the Energy Conservation of Polynomial Hamiltonian Problems.  {\em  JNAIAM J. Numer. Anal. Ind. Appl. Math.} {\bf 4} (2009) 87--101.

\bibitem{LeRe2004} B.\,Leimkuhler, S.\,Reich. {\em Simulating Hamiltonian Dynamics}. Cambridge University Press, Cambridge, UK, 2004.

\bibitem{LW2016}  Y.-W.\,Li, X.\,Wu. Functionally fitted energy-preserving methods for solving oscillatory nonlinear Hamiltonian systems. {\em SIAM J. Numer. Anal.} {\bf 54} (2016) 2036--2059.
\url{https://doi.org/10.1137/15M1032752}

\bibitem{LQR1999}  R.I.\,McLachlan, G.R.W.\,Quispel,  N.\,Robidoux. Geometric integration using discrete gradients. {\em Phil. Trans. R. Soc. Lond. A}  {\bf 357} (1999) 1021--1045.

\bibitem{SS2016} J.M.\,Sanz-Serna. Symplectic Runge-Kutta schemes for adjoint equations, automatic differentiation, optimal control, and more. {\em SIAM Rev.} {\bf 58} (2016) 3--33. \url{https://doi.org/10.1137/151002769}


\bibitem{SSC1994} J.M.\,Sanz-Serna, M.P.\,Calvo. {\em Numerical Hamiltonian Problems}. Chapman \& Hall, London, UK, 1994.

\bibitem{mrSIR} \url{https://www.mrsir.it/en/about-us/}

\end{thebibliography}
\end{document}